\documentclass[nospthms,envcountsect]{svmult} 
\usepackage{latexsym}
\usepackage{amsfonts}
\usepackage{amsmath}
\usepackage{amssymb}
\usepackage{amscd}
\usepackage{multicol}    
\usepackage{enumerate}
    \newtheorem{theorem}                    {Theorem}       [section]
    \newtheorem{lemma}      [theorem]       {Lemma}
    \newtheorem{corollary}  [theorem]       {Corollary}
    \newtheorem{proposition}[theorem]       {Proposition}
    \newtheorem{definition} [theorem]       {Definition}
    \newtheorem{conjecture} [theorem]       {Conjecture}

\begin{document}
\catcode`@=11
\atdef@ I#1I#2I{\CD@check{I..I..I}{\llap{$\m@th\vcenter{\hbox
  {$\scriptstyle#1$}}$}
  \rlap{$\m@th\vcenter{\hbox{$\scriptstyle#2$}}$}&&}}
\atdef@ E#1E#2E{\ampersand@
  \ifCD@ \global\bigaw@\minCDarrowwidth \else \global\bigaw@\minaw@ \fi
  \setboxz@h{$\m@th\scriptstyle\;\;{#1}\;$}%
  \ifdim\wdz@>\bigaw@ \global\bigaw@\wdz@ \fi
  \@ifnotempty{#2}{\setbox@ne\hbox{$\m@th\scriptstyle\;\;{#2}\;$}%
    \ifdim\wd@ne>\bigaw@ \global\bigaw@\wd@ne \fi}%
  \ifCD@\enskip\fi
    \mathrel{\mathop{\hbox to\bigaw@{}}%
      \limits^{#1}\@ifnotempty{#2}{_{#2}}}%
  \ifCD@\enskip\fi \ampersand@}
\catcode`@=\active

\renewcommand{\labelenumi}{\alph{enumi})}
\newcommand{\chr}{\operatorname{char}}
\newcommand{\isom}{\stackrel{\sim}{\longrightarrow}}
\newcommand{\Aut}{\operatorname{Aut}}
\newcommand{\Hom}{\operatorname{Hom}}
\newcommand{\End}{\operatorname{End}}
\newcommand{\HOM}{\operatorname{{\mathcal H{\frak{om}}}}}
\newcommand{\EXT}{\operatorname{\mathcal E{\frak xt}}}
\newcommand{\Tot}{\operatorname{Tot}}
\newcommand{\Ext}{\operatorname{Ext}}
\newcommand{\Gal}{\operatorname{Gal}}
\newcommand{\Pic}{\operatorname{Pic}}
\newcommand{\Spec}{\operatorname{Spec}}
\newcommand{\trdeg}{\operatorname{trdeg}}
\newcommand{\im}{\operatorname{im}}
\newcommand{\coim}{\operatorname{coim}}
\newcommand{\coker}{\operatorname{coker}}
\newcommand{\gr}{\operatorname{gr}}
\newcommand{\id}{\operatorname{id}}
\newcommand{\Br}{\operatorname{Br}}
\newcommand{\cd}{\operatorname{cd}}
\newcommand{\CH}{CH}
\newcommand{\Alb}{\operatorname{Alb}}
\renewcommand{\lim}{\operatornamewithlimits{lim}}
\newcommand{\colim}{\operatornamewithlimits{colim}}
\newcommand{\rk}{\operatorname{rank}}
\newcommand{\codim}{\operatorname{codim}}
\newcommand{\NS}{\operatorname{NS}}
\newcommand{\mot}{\operatorname{DM^{eff}_-}}
\newcommand{\cone}{{\rm cone}}
\newcommand{\rank}{\operatorname{rank}}
\newcommand{\ord}{{\rm ord}}
\newcommand{\f}{{\cal G}}
\newcommand{\g}{{\cal F}}
\newcommand{\du}{{\cal D}}
\newcommand{\G}{{\mathbb G}}
\newcommand{\N}{{\mathbb N}}
\newcommand{\A}{{\mathbb A}}
\newcommand{\Z}{{{\mathbb Z}}}
\newcommand{\Q}{{{\mathbb Q}}}
\newcommand{\R}{{{\mathbb R}}}
\newcommand{\B}{{\mathbb Z}^c}
\renewcommand{\H}{{{\mathbb H}}}
\renewcommand{\P}{{{\mathbb P}}}
\newcommand{\F}{{{\mathbb F}}}
\newcommand{\m}{{\mathfrak m}}
\newcommand{\Sch}{{\text{\rm Sch}}}
\newcommand{\et}{{\text{\rm et}}}
\newcommand{\Zar}{{\text{\rm Zar}}}
\newcommand{\Nis}{{\mathcal M}}
\newcommand{\tr}{\operatorname{tr}}
\newcommand{\tor}{{\text{\rm tor}}}
\newcommand{\red}{{\text{\rm red}}}
\newcommand{\Div}{\operatorname{Div}}
\newcommand{\Ab}{{\text{\rm Ab}}}
\newcommand{\DD}{{\mathbb Z}^c}
\renewcommand{\div}{\operatorname{div}}
\newcommand{\corank}{\operatorname{corank}}
\renewcommand{\O}{{\cal O}}
\newcommand{\C}{{\mathbb C}}
\newcommand{\p}{{\mathfrak p}}
\newcommand{\proof}{\noindent{\it Proof. }}
\newcommand{\proofend}{\hfill $Q.E.D.$ \\}
\newcommand{\rem}{\noindent {\it Remark. }}
\newcommand{\example}{\noindent {\bf Example. }}
\newcommand{\ar}{{\text{\rm ar}}}
\newcommand{\del}{{\delta}}
\title*{Finite generation conjectures for cohomology over finite fields}
\author{Thomas Geisser\thanks{Supported in part by NSF grant No.0901021}}
\institute{University of Southern California}

\maketitle

\vskip-1.5cm

\begin{abstract}
We relate motivic cohomology, and Weil-etale cohomology, both
of which are finitely generated, by an intermediate cohomology 
theory, which we also conjecture to be finitely generated, and
examine the relationship of the three theories.
\end{abstract}

\section{Introduction}
Bass conjectured that for a regular scheme $X$ of finite type
over the integers, the higher algebraic $K$-groups are finitely
generated \cite{bassln342}. 
Via the localization sequence, this is equivalent to 
the finite generation of $K'_i(X)$ for all $X$ of finite type
over the integers. In view of the spectral sequence from 
higher Chow groups to $K'$-theory, it is a slightly stronger
statement to conjecture the finite generation of higher Chow groups 
$CH_n(X,i)$ for all $X$ of finite type over the integers.
If we restrict ourselves to $X$ of finite type over a finite field,
then this is equivalent to the finite generation of 
motivic cohomology $H^i_\Nis(X,\Z(n))$ for smooth $X$
\cite{voevodsky}. Under resolution of singularities, it implies
finite generation of motivic cohomology for all $X$ of finite
type over a finite field, as one sees 
from the blow-up long exact sequence. 

Still for $X$ of finite type over a finite field, LIchtenbaum's 
Weil-etale cohomology $H^i_W(X,\Z(n))$ is defined as the 
cohomology of $R\Gamma_GR\Gamma_\et(\bar X,\Z(n))$, where $G$
is the Weil-group of the finite field. Using ideas of Kahn, we showed in \cite{ichweilI} 
that the finite generation of $H^i_W(X,\Z(n))$ for smooth
and proper $X$ is equivalent to the strong form of Tate's and 
Beilinson's conjecture. However, the groups 
$H^i_W(X,\Z(n))$ are not finitely generated if $X$ is not
smooth or not proper, and need to be modified in this case,
see \cite{ichweilII}.

The purpose of this article is to relate the two finite generaton conjectures
with the help of an intermediate cohomology theory $H^i_F(X,\Z(n))$, 
which we call Frobenius cohomology. It is defined as the cohomology of
$R\Gamma_G\Z(n)(\bar X)$, where $\Z(n)$ is Bloch's cycle complex
shifted appropriately. There are natural maps 
$$H^i_\Nis(X,\Z(n))\stackrel{\alpha}{\to}
H^i_F(X,\Z(n))\stackrel{g}{\to} H^i_W(X,\Z(n)),$$
which can be studied separately. In fact they lie
in a diagram
\begin{equation}\label{less}\begin{CD}
@>>>H^i_\Nis(X,\Z(n)) @>\alpha >>H^i_F(X,\Z(n)) 
@>\beta>> H^{i-1}_K(X,\Z(n)) @>\gamma>>\\
@EEE @Vf VV  @Vg VV @VhVV  \\
@>>> H^i_\et(X,\Z(n)) @>>> H^i_W(X,\Z(n)) @>>> H^{i-1}_\Nis(X,\Q(n))@>>>
\end{CD}\end{equation}
with groups $H^i_K(X,\Z(n))=H^i(\Z(n)(\bar X)_G)$ called Kato cohomology. 
They are a generalization of the integral version of Kato 
homology \cite{kato} 
(which is the case $n=\dim X$) defined in \cite{ichkato}. 
We conjecture all groups in the upper row of
\eqref{less} to be finitely generated,
and believe they form interesting invariants of $X$.
For example, $H^{2n}_K(X,\Z(n))=CH^n(\bar X)_G$, and  under Parshin's conjecture,
the groups $H^i_K(X,\Z(n))$ vanish for $i\leq n$ and are torsion
for $i\not= 2n$.
The first interesting example $H^3_K(X,\Z(2))$ consists of those
elements in the cokernel of the integral cycle map
$CH^2(X)\to H^4_\et(X,\Z(2))$ which are in the image over
the algebraic closure, or in classical language those elements
in the unramified cohomology group $H^3_{nr}(X,\Q/\Z(2))$ which vanish in 
$H^3_{nr}(\bar X,\Q/\Z(2))$. 
Pirutka \cite{pirutka} used an idea of 
Colliot-Th\'el\`ene to construct an element in $H^3_{nr}(X,\Z/2)$
for a geometrically rational surface (in all but finitely many
characteristics),  
providing a $2$-torsion element in $H^3_K(X,\Z(2))$.

Since the above theories exist for any scheme of finite
type over a finite field, we define them in this generality 
and later specialize to smooth projective $X$.

\section{Borel Moore homology theory}
We consider separated schemes of finite type over a 
finite field $\F_q$ of characteristic $p$, and $n\geq 0$
(the case $n<0$ can be reduced to this case using the homotopy formula).
Let $\varphi$ be the geometric Frobenius, and $G=\langle \varphi\rangle$
be the Weil group of $\F_q$.
Let $\Z^c(n)$ be the cycle complex defined by Bloch shifted by $2n$,
so that 
$$H_i(\Z^c(n)(X)) = CH_n(X,i-2n).$$
This agrees with motivic cohomology 
$H^{2d-i}_\Nis(X,\Z(d-n))$ for smooth $X$ of pure dimension $d$
by \cite{voevodsky}.

\begin{definition}\label{def1}
Let $A$ be an abelian group.
\begin{enumerate}
\item We define $H_i^c(X,A(n))$ to be $H_i(\Z^c(n)(X)\otimes A)$.
\item We define Frobenius homology $H_i^F(X,A(n))$ 
to be the homology of the double complex 
$$ \Z^c(n)(\bar X)\otimes A \xrightarrow{\varphi-1} \Z^c(n)(\bar X)\otimes A.$$
Here $\varphi$ acts covariantly in
$\bar X$ on cycles on $\bar X\times \Delta^i$; the left and right hand 
complexes sit in homological degrees $1$ and $0$, respectively.
\item The Kato homology $H_i^K(X,A(n))$ 
is defined to be the homology of the complex of coinvariants
$\big(\Z^c(n)(\bar X)\otimes A\big)_G$.
\end{enumerate}
\end{definition}

Note that $H_i^F(X,\Z/m(n))$ is isomorphic to 
$H^{1-i}(Gal(\F_q),\Z^c/m(n)(\bar X))$ because with torsion coefficients,
Galois cohomology can be calculated by the complex above.
The following lemma follows from the definitions.

\begin{lemma}\label{lem1}
The groups  $H_i^c(X,A(n))$, $H_i^F(X,A(n))$ and $H_i^K(X,A(n))$ 
vanish for $i<2n$.
We have $H_{2n}^c(X,A(n))\cong CH_n(X)\otimes A$ and
$$H_{2n}^K(X,A(n))\cong H_{2n}^F(X,A(n))
\cong \big( CH_n(\bar X)\otimes A\big)_G.$$
For all $i$, there are short exact sequences
\begin{equation}\label{yyy}
0\to H_i^c(\bar X,A(n))_G\to 
H_i^F(X,A(n))\to H_{i-1}^c(\bar X,A(n))^G \to 0.
\end{equation}
\end{lemma}

\begin{proposition}\label{prop1}
We have an exact sequence 
\begin{equation}\label{homd}
\cdots\to H_i^c(X,A(n)) \to H_{i+1}^F(X,A(n))\to H_{i+1}^K(X,A(n)) \to 
H_{i-1}^c(X,A(n)) \to \cdots 
\end{equation}
All three theories are covariently functorial for proper maps,
contravariantly functorial for quasi-finite flat maps,
and have localization long exact sequences.
\end{proposition}

\proof 
This comes from the short exact sequence of double complexes
$$\begin{CD}
\big( \Z^c(n)(\bar X)\otimes A\big)^G@>>> 
\Z^c(n)(\bar X)\otimes A@>>> 0 \\
@VVV @V\varphi-1 VV@VVV \\
0@>>> \Z^c(n)(\bar X)\otimes A @>>> 
\big(\Z^c(n)(\bar X)\otimes A\big)_G
\end{CD}$$
and $\Z^c(n)(X)\otimes A\cong \big(\Z^c(n)(\bar X)\otimes A\big)^G$.
Functoriality and the localization property are well-known
for cycle complexes, and this carries over to 
$\big( \Z^c(n)(\bar X)\otimes A\big)_G$ by an easy diagram chase.
\proofend

For a finitely generated field $K$ of transcendence degree $d$ over $\F_q$, 
define $H^i_F(K,A(j))$ to be $\colim H_{2d+1-i}^F(U,A(d-j))$,
where $U$ runs through smooth schemes with function field $K$.
By \eqref{yyy}, the cohomology groups of $K$ lie in short exact sequences
\begin{multline}\label{frobdescent}
0\to H^{s-t-1}_\Nis(K\otimes_{\F_q}\bar \F_q,A(s-n))_G
\to H^{s-t}_F(K,A(s-n))\to \\ 
H^{s-t}_\Nis(K\otimes_{\F_q}\bar \F_q,A(s-n))^G\to  0, 
\end{multline}
where $K\otimes_{\F_q}\bar \F_q$ is a finite product of fields.

\begin{proposition}\label{ssprop}
We have functorial spectral sequences
\begin{align*}
E^1_{s,t}=\bigoplus_{X_{(s)}}
H^{s-t}_\Nis(k(x),A(s-n))&\Rightarrow H_{s+t}^c(X,A(n));\\
\bar E^1_{s,t}=\bigoplus_{X_{(s)}}
H^{s-t}_F(k(x),A(s-n))&\Rightarrow H_{s+t+1}^F(X,A(n)).\\
\end{align*}
\end{proposition}

\proof 
This follows with the niveau filtration, together with 
the fact that 
$$H_{s+t}^c(K,A(n)):= \colim_{U} H_{s+t}^c(U,A(n))
\cong H^{s-t}_\Nis(K,A(s-n))$$
and $H_{s+t+1}^F(K,A(n))\cong H^{s-t}_F(K,A(s-n))$
for a field $K$ of transcendence degree $s$ over the base field.
\proofend

\subsection{Integral coefficients}
\begin{conjecture}\label{mainhom}
The groups $H_i^c(X,\Z(n))$, $H_{i}^F(X,\Z(n))$ and $H_{i}^K(X,\Z(n))$
are finitely generated for any $i, n$ and $X$.
\end{conjecture}

By localization, it suffices to consider the case of smooth and proper
$X$ (assuming that every finitely generated field over $\F_q$ 
has a smooth and proper model).
This case will be considered in detail below.

Recall from \cite{ichparshin} Parshin's conjecture $P_n$, stating that
$CH_n(X,i)$ is torsion for $i\not=0$ and $X$ smooth and projective.

\begin{proposition}\label{pupl}
Assume Conjecture $P_n$. Then
$$H^{s-t}_\Nis(k,\Z(s-n))\cong  H^{s-t}_F(k,\Z(s-n))$$
for $t\geq n$. For $t=n-1$, the left hand side vanishes, and 
$H^{s-n+1}_F(k,\Z(s-n))\cong K_{s-n}^M(k\otimes_{\F_q}\bar \F_q)_G$.
For $t<n-1$, both sides vanish.
\end{proposition}

\proof
It suffices to show that the second map and the composition
$$H^{s-t}_\Nis(k,\Z(s-n))\to H^{s-t}_F(k,\Z(s-n))
\to H^{s-t}_\et(k,\Z(s-n))$$
are isomorphism for $t\geq n$. 
The total composition is an isomorphism for $t\geq n-1$
by the Beilinson-Lichtenbaum conjectures. 
For the second map, we compare \eqref{frobdescent} to the
analog short exact sequence for etale cohomology. 
By Parshin's conjecture,
the groups $H^i_\Nis(k,\Z(s-n))$ are torsion for $i<s-n$,
so that Galois cohomology and cohomology of the Weil-group $G$
agree for those groups. The result now follows with
the Beilinson-Lichtenbaum conjectures by comparing \eqref{frobdescent}
with the Hochschild-Serre spectral sequence for Galois cohomology
(for the finite product of fields $k\otimes_{\F_q}\bar\F_q$).
\proofend

\begin{definition}
We define the Kato complex $KC_n(X)$ of weight $n$ to be the complex 
$$ \bigoplus_{X_{(n)}} \Z \gets 
\bigoplus_{X_{(n+1)}} (k(x)\otimes_{\F_q}\bar\F_q)^\times_G \gets
\cdots\gets 
\bigoplus_{X_{(s)}} K_{s-n}^M(k(x)\otimes_{\F_q}\bar\F_q)_G
\gets\cdots,$$
with the summand indexed by $X_{(s)}$ in degree $s-n$. The differentials
are induced by boundary maps of higher Chow groups of discrete
valuation rings.
\end{definition}

\begin{corollary}
Assuming conjecture $P_n$, we have 
$$H_i^K(X,\Z(n))\cong H_{i-2n}(KC_n(X)).$$
In particular, $H_i^K(X,\Z(n))$ vanishes unless $2n\leq i\leq n+d$.
\end{corollary}

\proof Compare the spectral sequences of Proposition \ref{ssprop}
and use Proposition \ref{pupl}.
\proofend

\subsection{Finite coefficients}
Considering finite coefficients allows us to remove
the hypothesis on Parshin's conjecture and to compare to 
etale homology.  Recall the sequence \eqref{homd}
\begin{equation}
\label{hommodm}
\cdots\to H_i^c(X,\Z/m(n)) \stackrel{\alpha}{\to} H_{i+1}^F(X,\Z/m(n)) 
\stackrel{\beta}{\to}  H_{i+1}^K(X,\Z/m(n))\to \cdots 
\end{equation}

Conjecture \ref{mainhom} implies the weaker

\begin{conjecture}\label{finitehom}
All terms in the \eqref{hommodm} are finite for all $i, n, X$.
\end{conjecture}

Let $f:X\to \F_q$ be the structure map.

\begin{definition}
For $p\not|m$ we define $H^\et_i(X,\Z/m(n))$ as the (Borel-Moore)
etale homology 
$H^{-i}(X_\et,Rf^!\mu_m^{\otimes -n})=
H^{-i}(Gal(\F_q),R\Gamma_\et(\bar X, Rf^!\mu_m^{\otimes -n}))$.
\end{definition}

Since $Gal(\F_q)$ has cohomological dimension $1$ and 
$Rf^!\mu_m^{\otimes -n}$ is concentrated in negative degrees, we have

\begin{lemma}
The groups $H^\et_i(X,\Z/m(n))$ vanish for $i<-1$.
\end{lemma}

Writing $\Z^c/m(n)$ for $\Z^c(n)\otimes \Z/m$, we showed in  
\cite{ichduality}:

\begin{theorem}\label{dduuaall}
We have $Rf^!\Z/m \cong \Z^c/m(0)$ for any $m$. In particular
$H^\et_i(X,\Z/m(n))$ agrees with the $-i$th etale hypercohomology of 
$\Z^c/m(n)$ if $p\not|m$ and $n\leq 0$.
\end{theorem}

If $p\not| m$ and $k$ contains the
$m$-th roots of unity, then cap-product with $\Z/m(1)(k)\cong \mu_m(k)$,  
defines a map $f^*\mu_m\otimes\Z^c/m(n)\to \Z^c/m(n-1)$
of complexes of sheaves on $X/k$, inducing a map
$$ \Z^c/m(n)\to f^*\mu^{\otimes -n}_m\otimes \Z^c/m(0)
\cong f^* \mu^{\otimes -n}_m\otimes Rf^!Z/m \cong Rf^!\mu^{\otimes -n}_m.$$
This map is in general not an isomorphism for $n>0$;  for example,
the left hand side vanishes for $n>\dim X$, but the right hand side
is periodic in $n$.
This construction induces a a map 
$$ \Z^c/m(n)(\bar X) \to R\Gamma_\et(\bar X, \Z^c/m(n))\to 
R\Gamma_\et(\bar X,Rf^!\mu^{\otimes -n}),$$
hence 
$$\gamma: H_{i+1}^F(X,\Z/m(n))\to H_i^\et(X,\Z/m(n))$$
by taking the cone of $\varphi-1$. 
The shift in degrees stems
from homological notation for Galois cohomology for the former, 
and cohomological notation for the latter. 

\begin{theorem}\label{KJS} Let $n=0$. Then 
Conjecture \ref{finitehom} holds for $p\not| m$, and for any $m$
under resolution of singularities.
\end{theorem}

\proof
If $n=0$, then $\gamma$ is an isomorphism by \cite{ichduality}, hence
Frobenius homology is finite. But by Jannsen-Kerz-Saito \cite{kerzsaito},
Kato homology is finite. 

If $m=p^r$, then Kato homology is finite under resolution of 
singularities. To show finiteness of Frobenius homology, one reduces by the usual devissage to the case that $X$ is smooth and proper of dimension $d$. 
In this case, $\gamma$ is an isomorphic to the finite group
$H^{d-i}_\et(X,\nu^{d}_r)$ by Theorem \ref{dduuaall} and the isomorphism 
$\Z^c/p^r(0) \cong \nu^d_r[d]$.
\proofend

As in Proposition \ref{ssprop}, we obtain 

\begin{proposition} 
For $p\not| m$ we have a spectral sequence
\begin{equation}\label{polk}
\tilde E^1_{s,t}=\bigoplus_{X_{(s)}}
H^{s-t}_\et(k(x),\mu_m^{\otimes s-n})
\Rightarrow H_{s+t}^\et(X,\Z/m(n))
\end{equation}
The $\tilde E^1$-terms lie in short exact sequences
\begin{multline}\label{galoismodm}
0\to H^{s-t-1}_\et(k(x)\otimes_{\F_q}\bar \F_q,\mu_m^{\otimes s-n})_G
\to H^{s-t}_\et(k(x),\mu_m^{\otimes s-n})\to \\ 
H^{s-t}_\et(k(x)\otimes_{\F_q}\bar \F_q,\mu_m^{\otimes s-n})^G\to  0. 
\end{multline}
\end{proposition}

\begin{proposition}\label{polk1}
For $k$ a field of transcendence degree $s$ over $\F_q$ and $p\not|m$, 
$$H^{s-t}_\Nis(k,\Z/m(s-n))\cong  H^{s-t}_F(k,\Z/m(s-n))
\cong H^{s-t}_\et(k,\mu_m^{\otimes s-n})
$$
for $t\geq n$. The left term vanishes for $t<n$, the middle
term vanishes for $t<n-1$ and 
$$H^{s-n+1}_F(k,\Z/m(s-n))\cong 
H^{s-n}_\Nis(k\otimes_{\F_q}\bar \F_q,\Z/m(s-n))_G
\subseteq H^{s-n+1}_\et(k,\mu_m^{\otimes s-n})$$
The right term vanishes for $t<-1$.
\end{proposition}

\proof
It suffices to show that the maps 
$$H^{s-t}_\Nis(k,\Z/m(s-n))\to H^{s-t}_F(k,\Z(s-n))
\to H^{s-t}_\et(k,\Z/m(s-n))$$
are isomorphism for $t\geq n$. The total composition is an isomorphism 
by the Beilinson-Lichtenbaum conjectures. Comparing the
sequences \eqref{frobdescent} and \eqref{galoismodm},
it follows that the second map
is an isomorphism by the Beilinson-Lichtenbaum conjectures
for the finite product of fields $k\otimes_{\F_q}\bar\F_q$.
\proofend

\begin{definition}
We define $KC_n/m(X)$ to be the complex 
\begin{multline}\label{kkk}
\bigoplus_{X_{(n)}} \Z/m \gets 
\bigoplus_{X_{(n+1)}} H^1_\et(k(x)\otimes_{\F_q}\bar\F_q,\Z/m(1))_G \gets
\cdots \\ \gets 
\bigoplus_{X_{(s)}} H^{s-n}_\et(k(x)\otimes_{\F_q}\bar\F_q,\Z/m(s-n))_G
\gets\cdots,
\end{multline}
with the summand indexed by $X_{(s)}$ in degree $s-n$.
\end{definition}

As in the integral case we obtain:

\begin{corollary}
We have 
$$H_i^K(X,\Z/m(n))\cong H_{i-2n}(KC_n/m(X)),$$
and theses groups vanish unless $2n\leq i\leq n+d$.
\end{corollary}

\section{Smooth and proper $X$}
In this section we assume that $X$ is smooth and proper
over $\F_q$, and discuss what other conjectures (like
Parshin's conjecture) imply for the cohomology groups
defined above. We use cohomological notation because
readers may be more familiar with it. 
Let $\Z(n)$ be Bloch's higher
cycle complex indexed by codimension, so that $\Z(n)(X)^i=z^n(X,2n-i)$
and $H^{i}_{Zar}(X,\Z(n))=CH_{d-n}(X,2n-i)$ is motivic cohomology for 
smooth $X$ of pure dimension $d$.
If $X$ is smooth of pure dimension $d$, we set
\begin{align}\label{definecohom}
H^i_F(X,A(n))&=H_{2d+1-i}^F(X,A(d-n))\\
H^i_K(X,A(n))&=H_{2d-i}^K(X,A(d-n)).
\end{align}

The following are reformulations of Lemma \ref{lem1}:

\begin{lemma}\label{plm}
The groups $H^i_F(X,A(n))$ vanish for $i>2n+1$, and in general
there are short exact sequences
$$ 0\to H^{i-1}_\Nis(\bar X,A(n))_G\to 
H^i_F(X,A(n))\to H^i_\Nis(\bar X,A(n))^G \to 0.$$
The groups $H^i_K(X,A(n))$ vanish for $i>2n$, and 
$$H^{2n}_K(X,A(n))\cong H^{2n+1}_F(X,A(n))
\cong \big(CH^n(\bar X)\otimes A\big)_G.$$
\end{lemma}

\begin{proposition}
We have a commutative diagram with exact rows
\begin{equation}\label{maind}\begin{CD}
@>>>{\bf H^i_\Nis(X,\Z(n))} @>\alpha >>{\bf  H^i_F(X,\Z(n))} 
@>\beta>> {\bf H^{i-1}_K(X,\Z(n))} @>\gamma>>\\
@EEE @Vf VV  @Vg VV @VhVV  \\
@>>> H^i_\et(X,\Z(n)) @>>> {\bf H^i_W(X,\Z(n)) }@>>> H^{i-1}_\Nis(X,\Q(n))@>>>
\end{CD}\end{equation}
\end{proposition}

\proof
The upper row is exact by Proposition \ref{prop1}, and the lower row
is the exact sequence relating etale cohomology to Weil-etale
cohomology \cite{ichweilI}. The maps $f$ and $g$ are induced by 
the change of topology from the Zariski to the etale site.
\proofend

\rem
The vertical maps are isomorphisms after tensoring with $\Q$,
because motivic cohomology and etale hypercohomology of the 
motivic complex agree with rational coefficients.

\begin{conjecture}\label{maincohom}
The bold face terms are finitely generated for all $i,n$ and $X$.
\end{conjecture}

Finite generation of 
$H^i_\Nis(X,\Z(n))$ is a generalization of Bass' conjecture, and finite 
generation of $H^i_W(X,\Z(n))$ 
is equivalent to Tate's and Beilinson's conjecture.
Lichtenbaum conjectured \cite{lichtenbaummotiv} also finite generation for 
the etale cohomology groups, except 
$H^{2n+2}_\et(X,\Z(n)) \cong  H^2(G_k,H^{2n}_\et(\bar X,\Z(n)))$:

\begin{conjecture}\label{lichtconj}
The group $H^i_\et(X,\Z(n))$ is finite for all $i\not=2n,2n+1$,
finitely generated for all $i=2n$, and of cofinite type for $i=2n+2$.
\end{conjecture}

\begin{lemma}\label{lowdegree}
The map $h:H^i_K(X,\Z(n))\to H^i_\Nis(X,\Q(n))$ 
is bijective for $i<n$ and injective for $i=n$.
\end{lemma}

\proof
By the Beilinson-Lichtenbaum conjectures, both the map 
$f$ as well as the map $g$ are isomorphism for $i\leq n+1$
and injective for $i=n+2$.
\proofend

\subsection{The case $n=0, 1, 2, d$}
Conjecture \ref{maincohom} holds for $n=0$. In fact we have
$H^i_F(X,\Z)\cong  H^i_W(X,\Z)\cong \Z$ for $i=0,1$,
$H^0_K(X,\Z)\cong H^0_\Nis(X,\Z)\cong \Z$, and all 
other homology groups vanish.

\begin{proposition} Let $n=1$. 

a) For $i=1$, the four left groups in \eqref{maind} are isomorphic to 
${\mathcal O}_X(X)^\times$, and for $i=2$, they are isomorphic
to $\Pic(X)$.

b) The groups $H^i_K(X,\Z(1))$ vanish for $i\not=2$, and
$H^2_K(X,\Z(1))\cong \NS(\bar X)_G$.

c) For $i=3$, 
\begin{align*}
H^3_\Nis(X,\Z(1))&=0\\
H^3_\et(X,\Z(1))&\cong \Br(X)\\
H^3_F(X,\Z(1))&\cong \NS(\bar X)_G,
\end{align*}
and $H^3_W(X,\Z(1))$ is an extension of $Br(\bar X)^G$
by $NS(\bar X)_G$. 

d) All groups in \eqref{maind} are finitely generated,
except possibly $H^3_\et(X,\Z(1))\cong \Br(X)$ and 
$H^3_W(X,\Z(1))$, 
whose finiteness is equivalent to Tate's conjecture for divisors.
\end{proposition} 

\proof This follows from $\Z(1)\cong {\mathbb G}_m[-1]$, Lemma \ref{plm}
and $(\Pic^0 \bar X)_G =0$.
\proofend

Now consider the case $n=2$. According to Parshin's conjecture
and Lemma \ref{lowdegree}, the groups $H^i_K(X,\Z(2))$ vanish for $i\leq 2$,
and $H^4_K(X,\Z(2))\cong CH^2(\bar X)_G$. To describe the remaining
group, let $\epsilon:X_{\et}\to X_{\Zar}$ be the canonical map of sites, and 
$\tau$ the composition of the boundary map in 
the lower row of \eqref{maind} with the composition
$$H^4_\et(X,\Z(2))\to H^0(X,R^4\epsilon_*\Z(2))\cong H^0(X,R^3\epsilon_*\Q/\Z(2))
=H^3_{nr}(X,\Q/\Z(2)).$$
The last term is known under the name unramified cohomology.

\begin{proposition}
The group ${}_\tor H^3_K(X,\Z(2))$ is isomorphic to the cohomology 
of the sequence
$$ H^2_\Nis(X,\Q(2))\stackrel{\tau}{\to} 
H^3_{nr}(X,\Q/\Z(2))\to H^3_{nr}(\bar X,\Q/\Z(2))^G.$$
\end{proposition}

\proof
The first statement follows from a diagram chase in \eqref{maind},
noting that due to the Beilinson-Lichtenbaum conjecture 
the maps $f$ and $g$ are injective,
$\coker f=H^0(X,R^4\epsilon_*\Z(2))$
and 
$$\coker g= \coker\big( H^4_\Nis(\bar X,\Z(2))^G\to 
H^4_\et(\bar X,\Z(2))^G\big) \subseteq H^0(\bar X,R^4\epsilon_*\Z(2))^G.$$
\proofend

\begin{corollary}
Under Parshin's conjecture there is an exact sequence 
$$0\to H^3_K(X,\Z(2))\to H^3_{nr}(X,\Q/\Z(2))\to 
H^3_{nr}(\bar X,\Q,\Z(2))^G .$$
Thus ${}_\tor H^3_K(X,\Z(2))$ consists of those elements
in the cokernel of the integral cyle map, whose pull back to 
the algebraic closure $H^4_\et(\bar X,\Z(2))$ lie in the image of
the cycle map $CH^n(\bar X)\to H^{2n}_\et(\bar X,\Z(n))$.
\end{corollary}

\proof
Under Parshin's conjecture, $H^3_k(X,\Z(2))$ is torsion. Furthermore, 
$\coker f$ is the obstruction to the integral
Tate conjecture, and $\coker g$ is the cokernel of the map
$CH^2(\bar X)^G\to H^4_\et(\bar X,\Z(2))^G$.
\proofend

\begin{corollary}
If $X$ is geometrically rational, then there are isomorphisms
$H^3_{nr}(X,\Q/\Z(2))\cong H^3_K(X,\Z(2))$.
\end{corollary}

\proof
For rational $Z$ over a field $k$, we have 
$H^3_{nr}(Z,\Q/\Z(2))=H^3_\et(k,\Q/\Z(2))$, and 
$H^2_\Nis(Z,\Q(2))\cong H^0(Z,{\mathcal H}^2(\Q(2)))\cong H^2_\Nis(k,\Q(2)) $
by \cite[2.1.9]{cothe}, and for $k$ the algebraic closure of a finite
field these groups vanish.
\proofend

Pirutka \cite{pirutka} 
constructed an element in $H^3_{nr}(X,\Z/2)$ for $X$ a geometrically
birational variety of dimension $5$, giving a non-trivial
$2$-torsion element in $H^3_K(X,\Z(2))$.

The following result is a consequence of the work of 
Jannsen, Kerz and Saito, see \cite{ichkato}:

\begin{theorem}
If $n=\dim X$, then Conjecture \ref{maincohom} is equivalent to 
conjecture $P_0$. In this case, $H^i_K(X,\Z(n))=0$ for $i\not=2n$,
and $H^{2n}_K(X,\Z(n))\cong \Z^{\pi_0(X)}$.
\end{theorem}

\subsection{Assuming Parshin's conjecture} 
The cohomological Parshin conjecture $P^n$ of \cite{parshinneu}
states that $H^i_\Nis(X,\Q(n))=0$ for smooth and proper $X$ and $i\not=2n$.
This is equivalent to isomorphisms 
$H^i_\et(X,\Z(n)) \cong H^i_W(X,\Z(n))$ for $i\leq 2n$, and injectivity
in degree $2n+1$ by \cite{ichweilI}. 
By Lemma \ref{lowdegree}, Parshin's conjecture implies that the groups 
$H^i_K(X,\Z(n))$ vanish for $i\leq n$ if $n>0$, and that they
torsion for $i\leq 2n$. 
Diagram \eqref{maind} for $i<2n$ becomes
\begin{equation}\label{ileq2n}
\begin{CD}
@>>> H^i_\Nis(X,\Z(n)) @>>> H^i_F(X,\Z(n)) 
@>>> H^{i-1}_K(X,\Z(n))\to \cdots \\
@III @VVV  @VVV   \\
@EEE H^i_\et(X,\Z(n)) @= H^i_W(X,\Z(n)) 
\end{CD}\end{equation} 
In degree $2n+1$ we get
$$\begin{CD}
0@>>> CH^n(\bar X)_\varphi @= H^{2n}_K(X,\Z(n)) @>>> 0 \\
@VVV @VVV  @VVV @VVV  \\
H^{2n+1}_\et(X,\Z(n))@>inj >> H^{2n+1}_W(X,\Z(n))@>>> CH^n(X)_\Q @>>>  
H^{2n+2}_\et(X,\Z(n))
\end{CD}$$  
The cokernel of the lower right horizontal map is the (conjecturally 
finite) group $H^{2n+2}_W(X,\Z(n))$. 
In degree $2n$ we have
\begin{equation}\label{degtwon}
\begin{CD}
H^{2n-2}_K(X,\Z(n))@>>> CH^n(X) @>\alpha >> H^{2n}_F(X,\Z(n)) 
@>>> H^{2n-1}_K(X,\Z(n))\to  0\\
 @III @Vf VV  @Vg VV   \\
@EEE H^{2n}_\et(X,\Z(n)) @= H^{2n}_W(X,\Z(n)) 
\end{CD}\end{equation}
We get an exact sequence
$$H^{2n-2}_K(X,\Z(n))\to \ker f\to \ker g\to H^{2n-1}_K(X,\Z(n))
\to \coker f\to \coker g \to 0.$$
By Corollary \ref{tatecor} below, $f$ has conjecturally the same kernel 
and cokernel as the cycle map, whereas $g$ is related to the cycle map
over the algebraic closure:
$$ \begin{CD}
0@>>> H^{2n-1}_\Nis(\bar X,\Z(n))_G @>>> 
H^{2n}_F(X,\Z(n))@>>> CH^n(\bar X)^G @>>> 0\\
@III @VVV @Vg VV @VVV \\
0@>>> H^{2n-1}_\et(\bar X,\Z(n))_G @>>> 
H^{2n}_W(X,\Z(n))@>>> H^{2n}_\et(\bar X,\Z(n))^G @>>> 0
\end{CD}$$
Thus Kato homology measures the difference of the failure of the integral
Tate conjecture over $\F_q$ and over its algebraic closure.






\section{The algebraic closure of a finite field}
As we saw, the map $H^i_\Nis(X,\Z(n))\to H^i_F(X,\Z(n))$
is controlled by $H^i_K(X,\Z(n))$. The maps 
$H^i_F(X,\Z(n))\to H^i_W(X,\Z(n))$
and $H^i_F(X,\Z/m(n))\to H^i_\et(X,\Z/m(n))$ arise
by taking Galois descent on the maps 
$$ H^i_\Nis(\bar X,\Z(n)) \to H^i_\et(\bar X,\Z(n)).$$
It is thus important to get some ideas of the properties of
this map. Since it is an isomorphism rationally, we focus on
finite coefficients.

Assume that $p\not|m$. By the proper base change theorem, 
$H^i_\et(\bar X,\Z/m(n))$ is finite, and replacing $\Z/m(n)$
by $\mu_m^{\otimes n}$, we have a localization long exact
sequence (the difference is that the latter is non-trivial for negative $n$). 
The question is if $H^i_\Nis(\bar X,\Z/m(n))$ is finite.
There are examples of Schoen showing that for certain threefolds
over an algebraically closed field of characteristic $0$, the group
$CH^2(\bar X)/l$ is not finite, which implies
that $H^4_\Nis(\bar X,\Z/l(2))$ cannot be finite. However, we are
not aware of any examples in characteristic $p$.

For $N\geq i$, consider the following diagram
$$\begin{CD}
H^i_\Nis(\bar X,\Z/m(n)) @>>> H^i_\et(\bar X,\mu_m^{\otimes n})\\
@V\cup \beta^{N-n} VV @| \\ 
H^i_\Nis(\bar X,\Z/m(N)) @= H^i_\et(\bar X,\mu_m^{\otimes N}).
\end{CD}$$
The lower row is an isomorphism by the Beilinson-Lichtenbaum
conjecture, and finiteness of $H^i_\Nis(\bar X,\Z/m(n))$ is 
equivalent to finiteness of the kernel of the cup-product with the Bott-element.

If $p=\chr k$, then $H^i_\et(\bar X,\Z/p(n))$ has no localization
long exact sequence.
The groups  $H^2_\et(\bar X,\Z/p(1))$ and 
$H^3_\et(\bar X,\Z/p(1))$ both contain the non-divisible $p$-torsion
of the Brauer group $H^3(\bar X_\et,\Z(1))$. 
For a supersingular abelian surface,
this group contains the field $\bar \F_q$, hence is not finitely
generated (however, its Galois invariants and coinvariants are
finite). According to Milne, the unipotent part of 
$H^i(X_\et,\Z/p(n))$ and $H^{2d+1-i}(X_\et,\Z/p(d-n))$ are in 
duality. 
There is a long exact sequence
$$\cdots  \to H^i_\Nis(\bar X,\Z/p(n)) \to H^i_\et(\bar X,\Z/p(n)) 
\to H^{i-1-n}_\Zar(\bar X,R^1\epsilon_* \nu^n)\to\cdots .$$
By the above example, $H^0_\Zar(\bar X,R^1\epsilon_* \nu^1)$
and $H^1_\Zar(\bar X,R^1\epsilon_* \nu^1)$ are infinite, but 
we don't know any example where $H^j_\Zar(\bar X,R^1\epsilon_* \nu^n)$ 
is infinite for $j\not= n-1, n$.

Regarding the integral structure of motivic cohomology, it is 
an interesting question if there is a presentation of the form
$$ H^i_{\mathcal M}(\bar X,\Z(n))\cong 
\Z^r\oplus (\Q/\Z')^c\oplus (\Q_p/\Z_p)^{c_p} \oplus (\text{finite})$$
with $r=0$ unless $i=2n$. In this case, 
$c$ would be independent of $n$ as soon as $n\geq i$, and it 
would be interesting to study the variation of $c$ for $n<i$.
The analog statement for etale cohomology is wrong
due to unipotent groups appearing, see the Brauer groups
above.

\section{The integral Tate conjecture}
In this section, $X$ is smooth and proper over a finite field $\F_q$.

\begin{proposition} Fix an integer $n$ and a scheme $X$, and 
consider the following statements:
\begin{enumerate} 
\item Lichtenbaum's conjecture \ref{lichtconj}.
\item The groups $H^i_W(X,\Z(n))$ are finitely generated.
\item Parshin's conjecture.
\end{enumerate}
Then $a)\Leftrightarrow b) \Rightarrow c)$.
\end{proposition}

\proof
$a)\Rightarrow c)$ This follows because
$H^i_\et(X,\Q(n))\cong H^i_\Nis(X,\Q(n))$, and the former vanishes
for $i\not=2n$ by hypothesis.


To show the equivalence of a) and b) we can assume Parshin's
conjecture, and consider the exact sequence of \cite{ichweilI}
$$ \cdots \to H^i_\et(X,\Z(n))\to H^i_W(X,\Z(n))\to H^{i-1}_\Nis(X,\Q(n))\to \cdots.$$
Then a) and b) imply each other for $i\leq 2n$
and $i>2n+2$, and we are left with 
\begin{multline*} 
0\to H^{2n+1}_\et(X,\Z(n))\to H^{2n+1}_W(X,\Z(n))\to H^{2n}_\Nis(X,\Q(n))\\
\to H^{2n+2}_\et(X,\Z(n))\to H^{2n+2}_W(X,\Z(n))\to 0 .
\end{multline*}
By the Weil-conjectures and counting coranks one sees that the corank $C$
of $H^{2n+2}_\et(X,\Z(n))\cong H^{2n+1}_\et(X,\Q/\Z(n))$ and of 
$H^{2n}_\et(X,\Q/\Z(n))$ agree.

a) $\Rightarrow$ b): Finiteness of $H^{2n+1}_\et(X,\Z(n))$ 
implies that $C$ agrees with the dimension of $H^{2n}_\et(X,\Q(n))$,
and since $H^{2n+2}_\et(X,\Z(n))$ is torsion,
$H^{2n+1}_W(X,\Z(n))/\tor$ must be a lattice of the same rank.
Moreover $H^{2n}_\Nis(X,\Q(n))$ surjects onto the divisible part
of $H^{2n+2}_\et(X,\Z(n))$, hence $H^{2n+2}_W(X,\Z(n))$ is finite.

b) $\Rightarrow$ a) We see that $H^{2n+1}_\et(X,\Z(n))$
is finite, because $H^{2n+1}_W(X,\Q(n))\cong H^{2n}(X,\Q(n))$
\cite{ichweilI}. For the same reason, $H^{2n+2}_\et(X,\Z(n))$
is an extension of the finite group $H^{2n+2}_W(X,\Z(n))$ and 
a torsion divisible group of finite corank.
\proofend

Let us relate the cycle map to the change of topology map for
motivic cohomology. Using the
functorial identification $\Z/m(n)\cong \mu_m^{\otimes n}$ 
for $\chr k\not| m$ of \cite{bkbl} and $\Z/p^r(n)\cong \nu_r^n[-n]$ of \cite{marcI}, 
the $l$-adic cycle map can be factored as follows:
$$ c:H^i_\Nis(X,\Z(n))\otimes\Z_l \stackrel{u}{\to} 
H^i_\et(X,\Z(n))\otimes\Z_l
\stackrel{v}{\to} H^i_\et(X,\Z(n))^{\wedge l} \stackrel{w}{\hookrightarrow} 
H^i_\et(X,\hat\Z_l(n)).$$
Here we write $A^{\wedge l}$ for the $l$-adic completion, and 
$H^i_\et(X,\hat\Z_l(n))$ for the $l$-adic
cohomology $\lim H^i_\et(X,\Z/l^r(n))$ to distinguish it from
the hypercohomology of $\Z(n)\otimes\Z_l$.

\begin{lemma}
The completion map $v$ is surjective.
\end{lemma}

\proof
The $\Z_l$-module  
$H^i_\et(X,\Z(n))^{\wedge l}\subseteq H^i_\et(X,\hat\Z_l(n))$
is finitely generated and the cokernel of $v$ is 
$l$-divsible by \cite[(4.2)]{jannsencont}, hence must be trivial.
\proofend

\begin{proposition}
We have an exact sequence
\begin{multline*}
0\to \ker u\to \ker c\to {\rm div}\; H^i_\et(X,\Z(n))\otimes\Z_l \\
\to \coker u\to \coker c\to T_lH^{i+1}(X_\et,\Z(n))\to 0.
\end{multline*}
The groups $\ker u$ and $\coker u$ are torsion.
\end{proposition}

\proof
This is the kernel-cokernel sequence of the composition
$c= (wv)\circ u$. The kernel of $wv$ is the kernel of $v$,
which is the group of divisible elements of $H^i_\et(X,\Z(n))\otimes \Z_l$.
The cokernel of $wv$ is the Tate module $T_lH^{i+1}(X,\Z(n))$.
\proofend

\begin{corollary}\label{tatecor}
Under Conjecture \ref{lichtconj}, the kernel and cokernel of $u$ are 
equal to the kernel and cokernel of the cycle map, respectively.
\end{corollary}

\proof
If $i\not=2n$, then all groups in the factorization of $c$
are expected to be finite, and $v$ and $w$ are isomorphisms.
The cycle map $CH^n(X)\otimes \Z_l\to H^{2n}_\et(X,\hat\Z_l(n))$
is rationally surjective if and only if 
$T_lH^{i+1}(X,\Z(n))$ vanishes, hence
Tate's conjecture is equivalent to $\Div H^{2n+1}_\et(X,\Z(n))$ 
torsion free.
It is rationally injective (Beilinson's
conjecture) if and only if $\div H^{2n}_\et(X,\Z(n))$ is torsion.
\proofend

Considering the diagram 
$$\begin{CD}
H^{2n-1}_\Nis(X,\Q/\Z(n))@>inj>> H^{2n}_\Nis(X,\Z(n))@>>> H^{2n}_\Nis(X,\Q(n))
@>surj>> H^{2n}_\Nis(X,\Q/\Z(n))\\
@VrVV @VuVV @| @VsVV \\
H^{2n-1}_\et(X,\Q/\Z(n))@>inj>> H^{2n}_\et(X,\Z(n))@>>> H^{2n}_\et(X,\Q(n))
@>>> H^{2n}_\et(X,\Q/\Z(n)),
\end{CD}$$
we obtain a short exact sequence
$0\to \coker r\to \coker u\to \ker s\to 0$, and all three groups
are conjecturally finite.

\subsection{The algebraic closure of a finite field}
Tate's original conjecture involved cycles over
the algebraic closure and taking invariants under the Galois group.
Let $X_r=X\times_{\F_q}\F_{q^r}$, $G_r=Gal(\bar \F_q/\F_{q^r})$
and consider the colimit over the composition
\begin{multline}\label{colimcycbar}
H^i_\Nis(X_r,\Z(n))\otimes\Z_l \stackrel{u}{\to} H^i_\et(X_r,\Z(n))\otimes\Z_l
\stackrel{v}{\twoheadrightarrow} H^i_\et(X_r,\Z(n))^{\wedge l}
\stackrel{w}{\hookrightarrow}\\ 
H^i_\et(X_r,\hat\Z_l(n))\stackrel{x}{\twoheadrightarrow}
H^i_\et(\bar X,\hat\Z_l(n))^{G_r}.
\end{multline}

It is clear that the first two groups commute with the colimit in $r$, so
if we write $\gamma_*M=\colim M^{G_r}$ for the largest
continuous submodule of a Galois-module, then 
the colimit takes the form 
\begin{multline}\label{colimcyc}
\bar c: H^i_\Nis(\bar X,\Z(n))\otimes\Z_l \stackrel{\bar u}{\to} 
H^i_\et(\bar X,\Z(n))\otimes\Z_l
\stackrel{\bar v}{\twoheadrightarrow} \colim H^i_\et(X_r,\Z(n))^{\wedge l}
\stackrel{\bar w}{\hookrightarrow}\\ 
\colim H^i_\et(X_r,\hat\Z_l(n))\twoheadrightarrow
\gamma_* H^i_\et(\bar X,\hat\Z_l(n)).
\end{multline}

\begin{corollary}
Under Conjecture \ref{lichtconj}, $\coker \bar c=\coker \bar u$.
\end{corollary} 

Note that surjectivity of $v$ implies that
$\colim H^i_\et(X_r,\Z(n))^{\wedge l}$ surjects onto
$H^i_\et(\bar X,\Z(n))^{\wedge l}$,
but the example of $\Br(X)=H^3_\et(X,\Z(1))$ shows that this surjection
is not an isomorphism in general.


\begin{thebibliography}{99}
\bibitem{bassln342} {\sc H.Bass}, Some problems in classical algebraic 
$K$-theory. Algebraic $K$-theory, II: ``Classical'' algebraic $K$-theory and 
connections with arithmetic (Proc. Conf., Battelle Memorial Inst., 
Seattle, Wash., 1972), pp. 3--73. Lecture Notes in Math., Vol. 342, 
Springer, Berlin, 1973.

\bibitem{bloch} {\sc S.Bloch}, Algebraic cycles and higher $K$-theory.
Adv. in Math. 61 (1986), no. 3, 267--304.

\bibitem{cothe} {\sc J.L. Colliot-Th\'el\`ene}. K-theory and algebraic geometry: 
connections with quadratic forms and division algebras (Santa Barbara, CA, 1992), 1-64, 
Proc. Sympos. Pure Math., 58, Part 1, Amer. Math. Soc., Providence, RI, 1995.

\bibitem{ichparshin} {\sc T.Geisser}, Parshin's conjecture revisited. 
$K$-theory and noncommutative geometry, 413--425, 
EMS Ser. Congr. Rep., Eur. Math. Soc., Zurich, 2008.

\bibitem{ichweilI} {\sc T.Geisser}, Weil-etale cohomology over finite fields. 
Math. Ann. 330 (2004), no. 4, 665--692. 

\bibitem{ichweilII} {\sc T.Geisser}, Arithmetic cohomology over finite fields 
and special values of $\zeta$-functions. 
Duke Math. J. 133 (2006), no. 1, 27--57.

\bibitem{ichduality} {\sc T.Geisser}, Duality via cycle complexes,
Ann. of Math. (2) 172 (2010), no. 2, 1095-1126. 

\bibitem{ichkato} {\sc T.Geisser}, Arithmetic homology and an integral version of 
Kato's conjecture, J. Reine Angew. Math. 644 (2010), 1-22.

\bibitem{parshinneu} {\sc T.Geisser}, Parshin's conjecture revisited,
$K$-theory and noncommutative geometry, 413--425, EMS Ser. Congr. Rep., 
Eur. Math. Soc., Zurich, 2008.

\bibitem{marcI} {\sc T.Geisser, M.Levine}, The $p$-part of
$K$-theory of fields in characteristic $p$. Inv. Math. 139 (2000), 459--494.

\bibitem{bkbl} {\sc T.Geisser, M.Levine}, The Bloch-Kato conjecture and a
theorem of Suslin-Voevodsky. J. Reine Angew. Math. 530 (2001), 55--103.

\bibitem{jannsencont} {\sc U.Jannsen}, 
Continuous etale cohomology. Math. Ann. 280 (1988), no. 2, 207--245.

\bibitem{jannsen} {\sc U.Jannsen}, Hasse principles for higher-dimensional
fields, Preprint Universit\"at Regensburg 18/2004.

\bibitem{jannsensaito} {\sc U.Jannsen, S.Saito}, Kato homology
of arithmetic schemes, and higher class field theory over local
fields, Doc. Math. (2003), 479-538.

\bibitem{kato} {\sc K.Kato}, A Hasse principle for two-dimensional global
fields. With an appendix by Jean-Louis Colliot-Th\'el\`ene. J. Reine
Angew. Math. 366 (1986), 142--183.

\bibitem{kerzsaito} {\sc M.Kerz, S.Saito}, Cohomological Hasse principle and
motivic cohomology for arithmetic schemes, Preprint 2010.

\bibitem{lichtenbaummotiv} {\sc S.Lichtenbaum}, 
Values of zeta-functions at nonnegative integers. Number theory, 
Noordwijkerhout 1983 (Noordwijkerhout, 1983), 127--138, 
Lecture Notes in Math., 1068, Springer, Berlin, 1984.

\bibitem{pirutka} {\sc A.Pirutka}, Sur le groupe de Chow de codimension
deux des vari\'et\'es sur les corps finis, Preprint 2010.

\bibitem{schneider} {\sc P.Schneider}, \"Uber gewisse Galoiskohomologiegruppen,
Math. Z. 168 (1979), no. 2, 181--205.

\bibitem{suslinetale} {\sc A.Suslin}, Higher Chow groups and
\'etale cohomology.
Cycles, transfers, and motivic homology theories, 239--254, Ann. of
Math. Stud., 143, Princeton Univ. Press, Princeton, NJ, 2000.

\bibitem{voevodsky} {\sc V.Voevodsky}, Motivic cohomology groups are
isomorphic to higher Chow groups in any characteristic. Int. Math.
Res. Not. 2002, no. 7, 351--355.
\end{thebibliography}
\end{document}